\documentclass[11pt]{amsart}
\usepackage{amsmath,amsthm,epsfig,amssymb}
\title{Platonic solids in $\mathbb Z^3$}
\author{Eugen J. Ionascu and Andrei Markov }
\curraddr{Department of Mathematics\\ Columbus State University\\4225 University Avenue\\
Columbus, GA 31907\\  First author is an honorific member of the
Romanian Institute of Mathematics ``Simion Stoilow" }
\email{ionascu\_eugen@colstate.edu, and markov\_andrei@colstate.edu}

\date{October $6^{th}$, 2009}
\begin{document}
\def\sms{\small\scshape}
\baselineskip18pt
%%%%%%%%%%%%%%%%%%%%%%%%%%%%%%%%%%%%%%%%%%%%%%%%%%%%%%%%%%%%%%%%%%%%%%
%%%%%%%%%%%%%%%%%%%% OUR DEFINITIONS %%%%%%%%%%%%%%%%%%%%%%%%%%%%%%%%%
%%%%%%%%%%%%%%%%%%%%%%%%%%%%%%%%%%%%%%%%%%%%%%%%%%%%%%%%%%%%%%%%%%%%%%
\newtheorem{theorem}{\hspace{\parindent}
T{\scriptsize HEOREM}}[section]
\newtheorem{proposition}[theorem]
{\hspace{\parindent }P{\scriptsize ROPOSITION}}
\newtheorem{corollary}[theorem]
{\hspace{\parindent }C{\scriptsize OROLLARY}}
\newtheorem{lemma}[theorem]
{\hspace{\parindent }L{\scriptsize EMMA}}
\newtheorem{definition}[theorem]
{\hspace{\parindent }D{\scriptsize EFINITION}}
\newtheorem{problem}[theorem]
{\hspace{\parindent }P{\scriptsize ROBLEM}}
\newtheorem{conjecture}[theorem]
{\hspace{\parindent }C{\scriptsize ONJECTURE}}
\newtheorem{example}[theorem]
{\hspace{\parindent }E{\scriptsize XAMPLE}}
\newtheorem{remark}[theorem]
{\hspace{\parindent }R{\scriptsize EMARK}}
\renewcommand{\thetheorem}{\arabic{section}.\arabic{theorem}}
\renewcommand{\theenumi}{(\roman{enumi})}
\renewcommand{\labelenumi}{\theenumi}
\newcommand{\Q}{{\mathbb Q}}
\newcommand{\Z}{{\mathbb Z}}
\newcommand{\N}{{\mathbb N}}
\newcommand{\C}{{\mathbb C}}
\newcommand{\R}{{\mathbb R}}
\newcommand{\F}{{\mathbb F}}
\newcommand{\K}{{\mathbb K}}
\newcommand{\D}{{\mathbb D}}
\def\phi{\varphi}
\def\ra{\rightarrow}
\def\sd{\bigtriangledown}
\def\ac{\mathaccent94}
\def\wi{\sim}
\def\wt{\widetilde}
\def\bb#1{{\Bbb#1}}
\def\bs{\backslash}
\def\cal{\mathcal}
\def\ca#1{{\cal#1}}
\def\Bbb#1{\bf#1}
\def\blacksquare{{\ \vrule height7pt width7pt depth0pt}}
\def\bsq{\blacksquare}
\def\proof{\hspace{\parindent}{P{\scriptsize ROOF}}}
\def\pofthe{P{\scriptsize ROOF OF}
T{\scriptsize HEOREM}\  }
\def\pofle{\hspace{\parindent}P{\scriptsize ROOF OF}
L{\scriptsize EMMA}\  }
\def\pofcor{\hspace{\parindent}P{\scriptsize ROOF OF}
C{\scriptsize ROLLARY}\  }
\def\pofpro{\hspace{\parindent}P{\scriptsize ROOF OF}
P{\scriptsize ROPOSITION}\  }
\def\n{\noindent}
\def\wh{\widehat}
\def\eproof{$\hfill\bsq$\par}
\def\ds{\displaystyle}
\def\du{\overset{\text {\bf .}}{\cup}}
\def\Du{\overset{\text {\bf .}}{\bigcup}}
\def\b{$\blacklozenge$}

\def\eqtr{{\cal E}{\cal T}(\Z) }
\def\eproofi{\bsq}

%%%%%%%%%%%%%%%%%%%%%%%%%%%%%%%%%%%%%%%%%%%%%%%%%%%
%%%%%%%%%%%%%%%%%%%%%%%ABSTRACT%%%%%%%%%%%%%%%%%%%%%%%%%%%%%%%%%%%%
\begin{abstract} Extending previous results on a
characterization of all equilateral triangle in space having
vertices with integer coordinates (``in $\mathbb Z^3$"),  we look at
the problem of characterizing all regular polyhedra (Platonic
Solids) with the same property. To summarize, we show first that
there is no regular icosahedron/ dodecahedron in $\mathbb Z^3$. On
the other hand, there is a finite (6 or 12) class of regular
tetrahedra in $\mathbb Z^3$, associated naturally to each nontrivial
solution $(a,b,c,d)$ of the Diophantine equation $a^2+b^2+c^2=3d^2$
and for every nontrivial integer solution $(m,n,k)$ of the equation
$m^2-mn+n^2=k^2$. Every regular tetrahedron in $\mathbb Z^3$
belongs, up to an integer translation and/or rotation, to one of
these classes. We then show that each such tetrahedron can be
completed to a cube with integer coordinates. The study of regular
octahedra is reduced to the cube case via the duality between the
two. This work allows one to basically give a description the
orthogonal group $O(3,\mathbb Q)$ in terms of the seven integer
parameters satisfying the two relations mentioned above.
\end{abstract} \maketitle
%%%%%%%%%%%%%%%%%%%%%%%%%%%%%%%%%%%%%%%%%%%%%%%%%%%%%%%%%%%%%%%%%%%%%%%
\section{INTRODUCTION}%%%%%%%%%%%%%%%%%%%%%%%%%%%%%%%%%%%%%%%%%%%%%%%%%
The set of equilateral triangles in the three dimensional space with
integers coordinates for its vertices is very rich. If one counts
all of these triangles with the coordinates in
$\{0,1,2,3,...,100\}$, finds that there are $10,588,506,416$ of them
(\cite{OL}). A constructive way to find these triangles is described
in \cite{rceji} and \cite{eji}, and this method was implemented in
Maple (see \cite{eji2}). In \cite{ejirt} a characterization of
regular tetrahedra in $\mathbb Z^3$ is given. One may ask naturally
if other Platonic solids may happen to have integer coordinates.
Perhaps the simplest different examples are those of cubes in space
whose vertices have integer coordinates can be found by dilating the
so called unit cube by an integer factor: $(0,0,0)$, $(n,0,0)$,
$(0,n,0)$, $(0,0,n)$, $(n,n,0)$, $(n,0,n)$, $(0,n,n)$, and
$(n,n,n)$, $n\in \mathbb Z$. The figure below shows a less obvious
cube which has vertices: $(1,4,3)$, $(3,3,1)$, $(1,1,0)$,
$(-1,2,2)$, $(2,2,5)$, $(4,1,3)$, $(2,-1,2)$ and $(0,0,4)$.

\begin{center}\label{fig1}
$\underset{Figure\ 1 (a): \ Non-trivial\
cube}{\epsfig{file=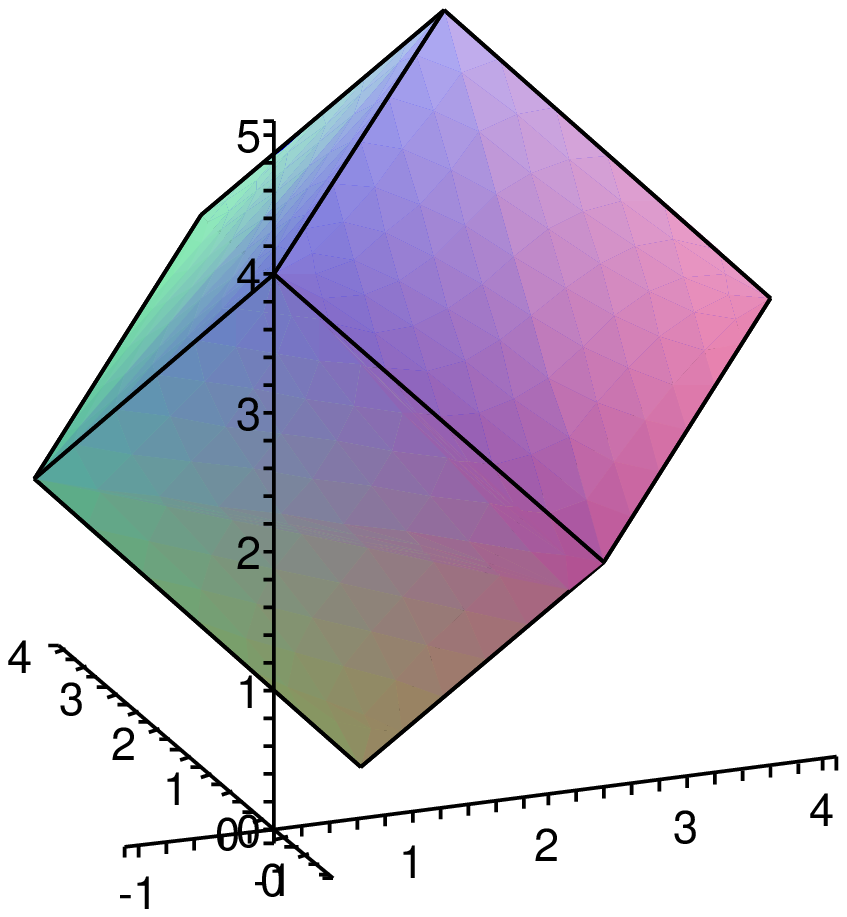,height=3in,width=3in}}\
\underset{Figure \ 1(b): \ Regular\ tetrahedron\ inscribed: \
OABC}{\epsfig{file=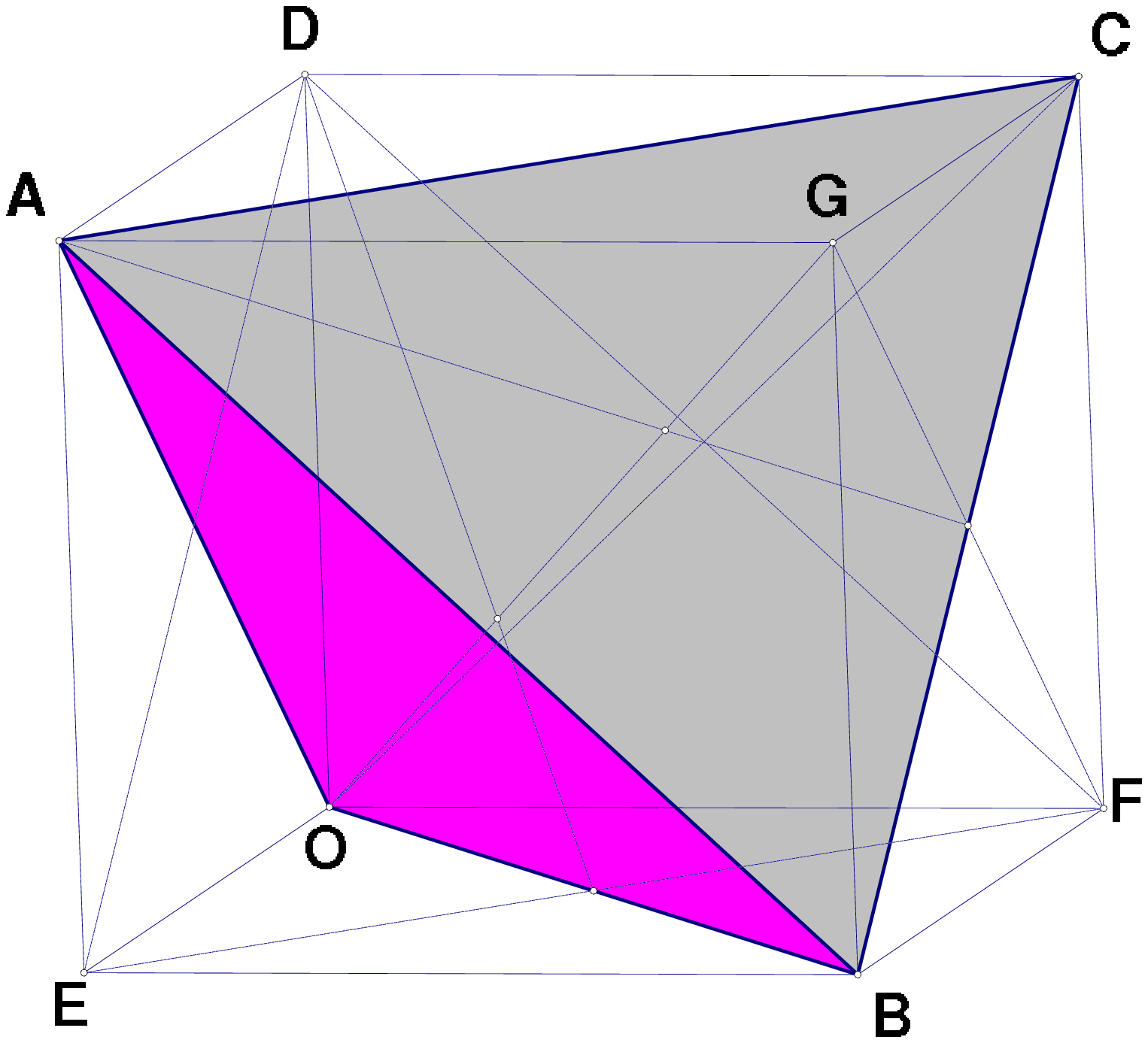,height=2in,width=2in}}$
\end{center}

On the other side of the question asked earlier, we will show that
dodecahedrons or icosahedrons do not exist ``in $\mathbb Z^3$".

 Given a cube in general, there is a natural way to imbed a
regular tetrahedron in it like in Figure 1(b). In fact, one can find
two such regular tetrahedra that have vertices a subset of the
cube's vertices but if we require that both the cube and the
tetrahedron have one of the vertices the origin, the correspondence
is one to one. If the cube has integer coordinates then the regular
tetrahedra considered must have the same property. Conversely, if we
start with a regular tetrahedra having integer coordinates this can
be completed uniquely by adding four more vertices, to a cube as in
Figure 1(b). In general, the resulting cube, as we will see, may
have coordinates only in $\frac{1}{2}\mathbb Z$, but it works out
that the coordinates are always in $\mathbb Z$.

Finally,  we study the regular octahedrons with integers coordinates
which are in duality with the cubes.

\section{Previous Results}
Our starting point is the theorem below about how to obtain all
equilateral triangles in ``in $\mathbb Z^3$". Each such triangle is
contained in a lattice of points of the form

\begin{equation}\label{planelattice}
{\cal P}_{a,b,c}:=\{(\alpha,\beta,\gamma)\in \mathbb {Z}^3|\ \
a\alpha+b\beta+c\gamma=0,\ \ a^2+b^2+c^2=3d^2,\ \ a,b,c,d\in \mathbb
Z\}.
\end{equation}

\begin{center}\label{figure0}
\epsfig{file=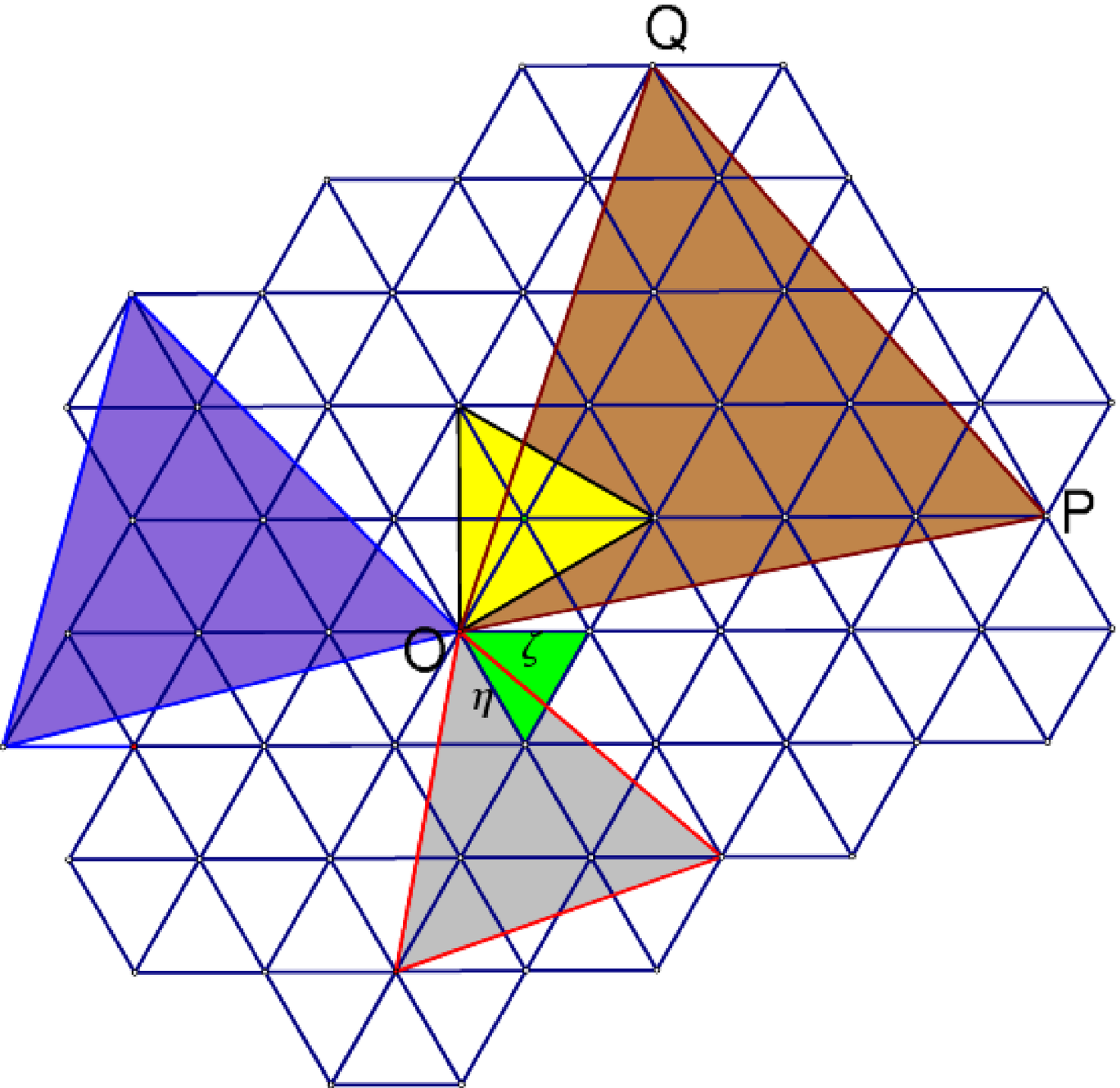,height=3in,width=3in}
\vspace{.1in}\\
{\small Figure 2: The lattice ${\cal P}_{a,b,c}$ } \vspace{.1in}
\end{center}

\n In general, the vertices of the equilateral triangles that dwell
in ${\cal P}_{a,b,c}$ form a strict sub-lattice of ${\cal
P}_{a,b,c}$  which is generated by only two vectors,
$\overrightarrow{\zeta}$ and $\overrightarrow{\eta}$ (see
Figure~\ref{figure0}). These two vectors are described by the
Theorem~\ref{generalpar} proved in \cite{rceji}.

\begin{theorem}\label{generalpar} Let $a$, $b$, $c$, $d$ be odd
integers such that $a^2+b^2+c^2=3d^2$ and ${\rm gcd}(a,b,c)=1$. Then
for every $m,n \in\mathbb{Z}$ (not both zero) the triangle $OPQ$,
determined by

\begin{equation}\label{vectorid}
\overrightarrow{OP}=m\overrightarrow{\zeta}+n\overrightarrow{\eta},\
\
\overrightarrow{OQ}=m(\overrightarrow{\zeta}-\overrightarrow{\eta})+n\overrightarrow{\zeta},
\ { \rm with} \ \overrightarrow{\zeta}=(\zeta_1,\zeta_1,\zeta_2),
\overrightarrow{\eta}=(\eta_1,\eta_2,\eta_3),
\end{equation}

\begin{equation}\label{paramtwo}
\begin{array}{l}
\begin{cases}
\ds \zeta_1=-\frac{rac+dbs}{q}, \\ \\
\ds \zeta_2=\frac{das-bcr}{q},\\ \\
\ds \zeta_3=r,
\end{cases}
,\ \
\begin{cases}
\ds \eta_1=-\frac{db(s-3r)+ac(r+s)}{2q},\\ \\
\ds \eta_2=\frac{da(s-3r)-bc(r+s)}{2q},\\ \\
\ds \eta_3=\frac{r+s}{2},
\end{cases}
\end{array}
\end{equation}
where $q=a^2+b^2$ and $(r,s)$ is a suitable solution of
$2q=s^2+3r^2$ that makes all the numbers in {\rm (\ref{paramtwo})}
integers, forms an equilateral triangle in $\mathbb Z^3$ contained
in the lattice {\rm (\ref{planelattice})} and having sides-lengths
equal to $d\sqrt{2(m^2-mn+n^2)}$.

Conversely, there exists a choice of the integers $r$ and $s$ such
that given an arbitrary equilateral triangle in $\mathbb{R}^3$ whose
vertices, one at the origin and the other two in the lattice  {\rm
(\ref{planelattice})}, then there also exist integers $m$ and $n$
such that the two vertices not at the origin are given by {\rm
(\ref{vectorid})} and {\rm (\ref{paramtwo})}.
\end{theorem}

Out of all the equilateral triangles in $\mathbb Z^3$ only those for
which $m^2-mn+n^2=\lambda ^2$, $\lambda \in \mathbb Z$, may give
rise to regular tetrahedra in $\mathbb Z^3$ according to the
following theorem. As a notation, we let for every $k\in \mathbb Z$,
$\Omega(k)=\{(m,n)\in \mathbb Z^2 |m^2-mn+n^2=k ^2\}$.

\begin{theorem}\label{main} Every tetrahedron whose side lengths are $\lambda\sqrt{2}$,
$\lambda\in \mathbb N$, which has a vertex at the origin,  can be
obtained by taking as one of its faces an equilateral triangle
having the origin as a vertex and the other two vertices given by
{\rm (\ref{vectorid})} and {\rm (\ref{paramtwo})} with $a$, $b$, $c$
and $d$ odd integers satisfying $a^2+b^2+c^2=3d^2$ with $d$ a
divisor of $\lambda$, and then completing it with the fourth vertex
$R$ with coordinates

\begin{equation}\label{fourthpoint}
\begin{array}{c}
 \displaystyle \left(
\frac{\begin{array}{c}
        (2\zeta_1-\eta_1)m \\
        -(\zeta_1+\eta_1)n \\
       \pm 2ak
      \end{array}
}{3},\frac{\begin{array}{c}
        (2\zeta_2-\eta_2)m \\
        -(\zeta_2+\eta_2)n \\
       \pm 2bk
      \end{array}}{3},
 \displaystyle \frac{\begin{array}{c}
        (2\zeta_3-\eta_3)m \\
        -(\zeta_3+\eta_3)n \\
       \pm 2ck
      \end{array}}{3} \right),  \ for\  some\  (m,n)\in
\Omega(k),\ k:=\frac{\lambda }{d}.
\end{array}
\end{equation}

\n Conversely, if we let $a$, $b$, $c$ and $d$ be a primitive
solution of $a^2+b^2+c^2=3d^2$, let $k\in \Bbb N$ and $(m,n)\in
\Omega(k)$, then the coordinates of the point $R$ in {\rm
(\ref{fourthpoint})}, which completes the equilateral triangle $OPQ$
given as in {\rm (\ref{vectorid})} and {\rm (\ref{paramtwo})}, are

(a) all integers, if $k\equiv 0$ {\rm (mod  3)} regardless of the
choice of signs or

(b) integers, precisely for only one choice of the signs if $k\not
\equiv 0$ {\rm (mod 3)}.
\end{theorem}

\section{New Facts}

As observed in \cite{ejirt} the dihedral angle between the faces of
a regular tetrahedron has a cosine of $\frac{1}{3}$. This is
compatible with the normals to equilateral triangles given
implicitly by (\ref{planelattice}). This is not the case for regular
icosahedra and regular dodecahedrons in $\mathbb Z^3$. The next
result may come to no surprise to our reader since there are so many
Diophantine equations that need to be satisfied.

\begin{center}\label{fig2}
$\underset{\small \ Figure\ 3\ (a): \ Regular\
Icosahedron}{\epsfig{file=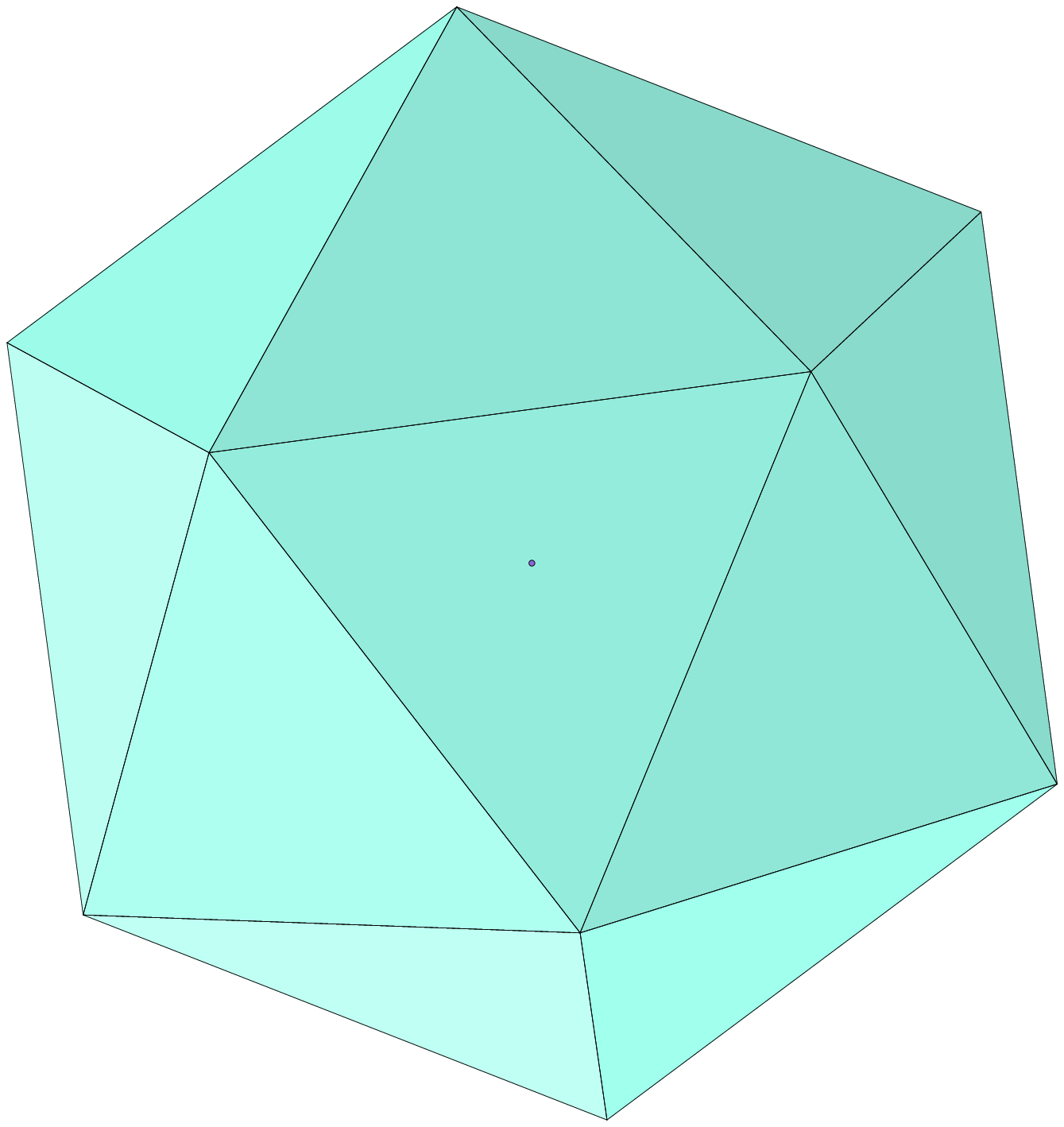,height=2in,width=2in}}$\
\ \  $\underset{\small \ Figure\ 3\ (b): \ Regular\
Dodecahedron}{\epsfig{file=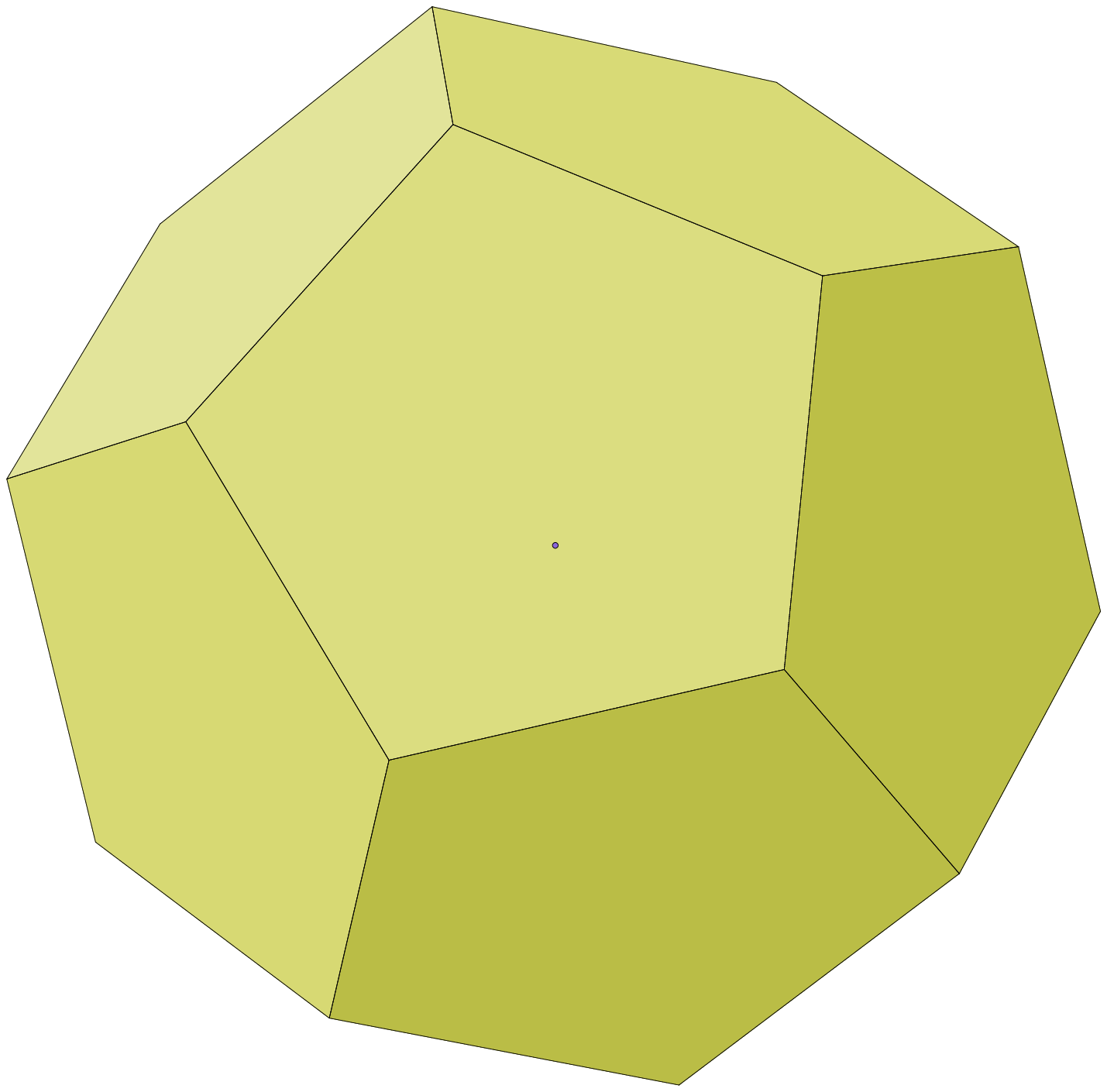,height=2in,width=2in}}$
\end{center}

\begin{theorem}\label{icosdodec} There is no regular icosahedron and no regular
dodecahedron in $\mathbb Z^3$.
\end{theorem}

\proof. \ We are beginning with the case of an icosahedron. We are
arguing by way of contradiction. It is public knowledge that the
cosine of the dihedral angle between the two adjacent faces of a
regular icosahedron is $-\frac{\sqrt{5}}{3}$ (which gives an angle
of approximately $138.1896851^{\circ}$). We consider two such
adjacent faces which are equilateral triangles (see Figure 3 (a))
and their normals: $\overrightarrow{n}=\frac{(a,b,c)}{d\sqrt{3}}$
and $\overrightarrow{n'}=\frac{(a',b',c')}{d'\sqrt{3}}$ for some odd
integers $a$, $b$, $c$, $d$, $a'$, $b'$, $c'$, and $d'$. But the
cosine of the two vectors above is equal to their scalar product
$\overrightarrow{n} \overrightarrow{n'}=\frac{aa'+bb'+cc'}{3dd'}\in
\mathbb Q$. This contradicts the fact that $\sqrt{5}$ is irrational.

For dodecahedrons, if there is one that is regular in $\mathbb Z^3$,
the dual polyhedron (the one obtained by taking the centers of mass
of each face), is a regular icosahedron which must have coordinates
that are rational numbers with denominators of 1 or 5. Translating
it so that a vertex becomes the origin and expanding everything by a
factor of five it becomes an regular icosahedron in $\mathbb Z^3$.
By what we have just shown this last object cannot exist. Hence,
there cannot be any regular dodecahedron in $\mathbb Z^3$. \eproof

Next let us use the Theorem~\ref{main} to give a characterization of
all cubes in $\mathbb Z^3$.

\begin{theorem}\label{cubes} Every cube in $\mathbb Z^3$ can be
obtained by a translation along a vector with integer coordinates
from a cube with a vertex the origin containing a regular
tetrahedron with a vertex at the origin and all integer coordinates
(see Figure 1(b)) and as a result it must have side lengths equal to
$n$ for some $n\in \mathbb N$. Conversely, given a regular
tetrahedron in $\mathbb Z^3$, this can be completed to a cube which
is going to be automatically in $\mathbb Z^3$.
\end{theorem}

\proof.  The first part of this statement follows from
Theorem~\ref{main} which prescribes the sides of a regular
tetrahedron to be of the form $n\sqrt{2}$, $n\in \mathbb N$. The
converse requires other parts of Theorem~\ref{main}. We may assume
that vertices of the given regular tetrahedron are given by the
origin and three other points, say $A(a_1,a_2,a_3)$,
$B(b_1,b_2,b_3)$ and $C(c_1,c_2,c_3)$. By Theorem~\ref{main}, we may
assume that $A$, $B$ are given by {\rm (\ref{vectorid})}, {\rm
(\ref{paramtwo})} and $C$ is given by (\ref{fourthpoint}) (Figure
2(b)). We observe that the planes defined by ABC and DEF are
parallel and cut the diagonal $\overline{OG}$  into three equal
parts. Hence the coordinates of $G$  must be
$(\frac{a_1+b_1+c_1}{2},\frac{a_2+b_2+c_2}{2},\frac{a_3+b_3+c_3}{2})$.
Similarly, the coordinates of $D$, $E$ and $F$ are respectively
$(\frac{a_1-b_1+c_1}{2},\frac{a_2-b_2+c_2}{2},\frac{a_3-b_3+c_3}{2})$,
\
$(\frac{a_1+b_1-c_1}{2},\frac{a_2+b_2-c_2}{2},\frac{a_3+b_3-c_3}{2})$,
and
$(\frac{-a_1+b_1+c_1}{2},\frac{-a_2+b_2+c_2}{2},\frac{-a_3+b_3+c_3}{2})$.
It suffices to show that the coordinates of $G$ are integers. Let us
look at $a_1+b_1+c_1$ (mod 2). We have
$$a_1+b_1+c_1=m\zeta_1+n\eta_1+m(\zeta_1-\eta_1)+n\zeta_1+(2\zeta_1-\eta_1)m-(\zeta_1+\eta_1)n\mp 2ak\equiv 0\  (mod\  2).$$
By symmetry the other two components must satisfy similar
relations.\eproof

We will refer to cubes or other Platonic solids in $\mathbb Z^3$ as
being {\it irreducible} if it cannot be obtained by an integer
vector translation and  multiplication  by an integer factor,
greater than one in absolute value,  from another analog object in
$\mathbb Z^3$.

\begin{corollary}\label{length} Every irreducible cube must have
side lengths which are only odd natural numbers.
\end{corollary}

\proof. \ If a given cube is irreducible, then it must arrive from
an irreducible tetrahedron $\cal T$. Indeed, if the formulae for the
vertices of a  tetrahedron $\cal T$ can be simplified to a smaller
tetrahedron $\cal T'$ then we can use the construction above and
obtain that the given cube is obtained from the cube corresponding
to $\cal T'$. Hence, we may assume that the tetrahedron $\cal T$ is
irreducible, and by Theorem~\ref{generalpar} and Theorem~\ref{main},
its side lengths are $d\sqrt{2(m^2-mn+n^2)}$ with $m$ and $n$
relative prime numbers such that $m^2-mn+n^2=k^2$, $k\in \mathbb N$
and $d$ an odd number. This implies that that the sides of the given
cube are equal to $dk$. Here, $k$ must be also odd, otherwise $m$
and $n$ must be divisible by $2$.\eproof

The problem of finding the number of cubes in space with coordinates
in $\{0,1,2,....,n\}$ has been studied in \cite{il}. In
\cite{ejict&c}, we implemented the method of obtaining cubes from
tetrahedra and the former from equilateral triangles. So we extended
the sequences A098928 and A103158. We list here a few more terms in
the sequence A098928.

\vspace{0.1in}

\centerline{ \vspace{0.2in}
\begin{tabular}{|c||c|c|c|c|c|c|c|c|c|c|c|}
  \hline
  n& 1 & 2 & 3 & 4 & 5 & 6 & 7 & 8 & 9 &10&11\\ \hline
 A098928& 1 & 9 & 36 & 100 & 229 & 473 & 910 & 1648 & 2795 & 4469 & 6818  \\
  \hline
\end{tabular}}
\vspace{0.1in}

\centerline{
\begin{tabular}{|c||c|c|c|c|c|c|c|}
  \hline
  n &12&13& 14& 15 & 16 & 17 & 18 \\ \hline
 A098928&10032 & 14315& 19907& 27190 & 36502 & 48233 & 62803 \\
  \hline
\end{tabular}.}
\vspace{0.1in}

It is clear that $A098928 \le A103158$. For $n\ge  4$ we have
actually a strict inequality, $A098928 < A103158$, and this  is due
to the fact that some of the tetrahedrons  inside of grid
$\{0,1,..,n\}^3$ extend beyond of its  boundaries.

\begin{center}\label{fig2}
$\underset{\small \ Figure\ 4\ (a) : \ Regular\
Octahedron}{\epsfig{file=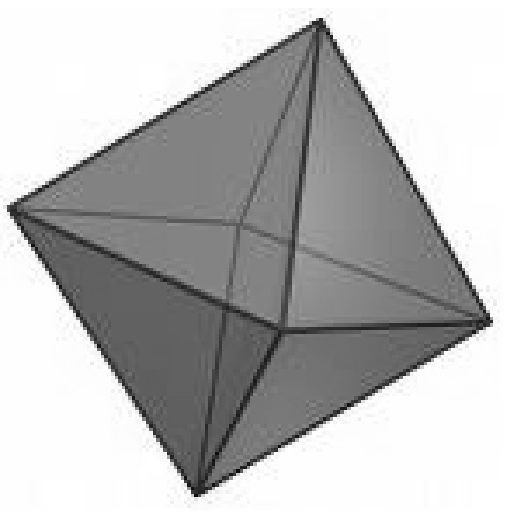,height=2in,width=2in}}$
$\underset{\small \ Figure\ 4\ (b) : \ Regular\
Octahedron}{\epsfig{file=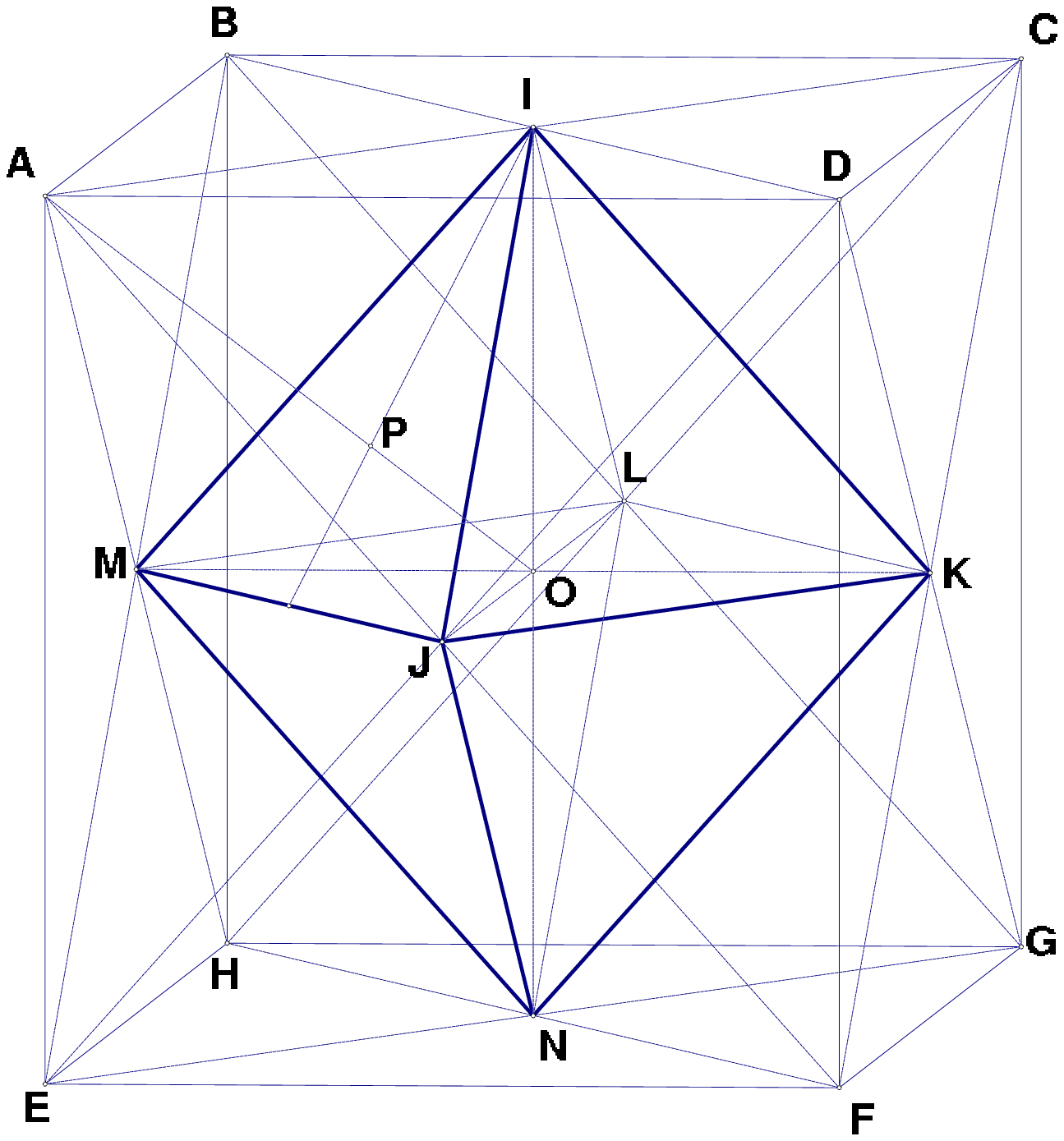,height=2in,width=2in}}$
\end{center}

Regular octahedrons in $\mathbb Z^3$ can be obtained from cubes in
$\mathbb Z^3$, doubling  their size and then taking the centers of
the faces. The converse of this statement is also true.

\begin{theorem} Every regular octahedrons in $\mathbb Z^3$ is the
dual of a cube that can be obtained (up to a translation with a
vector with integer coordinates) by doubling a cube in $\mathbb
Z^3$.
\end{theorem}

\proof. \ As we can see from Figure 4 (b), assuming that the
octahedron IJKLMN has vertices with integer coordinates, the centers
of the faces have coordinates that are rational numbers in
$\frac{1}{3}\mathbb Z$.   We may assume without loss of generality
that the center of the octahedron is the origin. In the Figure 4
(b), the octahedron IJKLMN is the dual of the cube ABCDEFGH. The
point $P$ is the center of the face IJM, and from similarity of
triangles we see that the coordinates are actually three times the
coordinates of $P$. This shows that the cube ABCDEFGH is in $\mathbb
Z^3$. Then the tetrahedron IJMA is then in $\mathbb Z^3$ and so by
Theorem~\ref{cubes} can be completed to a cube in $\mathbb Z^3$.
This cube if dilated by  a factor of two can be shifted to get the
cube ABCDEFGH.\eproof

If we denote by $\cal{RO}(n)$, the number of regular octahedron
whose vertices are in the set $\{0,1,...,n\}^3$, we get the
following sequence:

\centerline{ \vspace{0.2in}
\begin{tabular}{|c||c|c|c|c|c|c|c|c|c|c|c|}
  \hline
  n& 1 & 2 & 3 & 4 & 5 & 6 & 7 & 8 & 9 &10&11\\ \hline
$\cal{RO}(n)$ & 0 & 1 & 8 & 32 & 104 & 261 & 544 & 1000 & 1696 & 2759 & 4296  \\
  \hline
\end{tabular}}
\vspace{0.1in}

\centerline{
\begin{tabular}{|c||c|c|c|c|c|c|c|}
  \hline
  n &12&13& 14& 15 & 16 & 17 & 18 \\ \hline
$\cal{RO}(n)$ &6434 & 9352& 13243& 18304 & 24774 & 32960 & 43223 \\
  \hline
\end{tabular}.}

\vspace{0.1in}

These constructions bring up the idea of a certain order on the
orthogonal matrices with rational entries. For example, the first
level of such matrices are the ones in which the entries are
actually integers. The second tier will be formed by  those matrices
with denominators equal to 3, such as

$$T _3:= \frac{1}{3}
\left( \begin{array}{ccc}
          1 & -2 & 2 \\
          2 & -1 & -2 \\
          -2 & -2 & -1 \\
        \end{array}
      \right).
$$

Taking the square of $T_3$ we get a matrix in tier nine. There are
only two representations of the type
$3(9^2)=5^2+7^2+13^2=1^2+11^2+11^2$, and these give the same
``class" of matrices in the tier nine.  The first case when we have
essentially two different classes is for tier 13. An example of an
orthogonal matrix with denominators equal to 2009 is given below:

$$T_{2009} := \frac{1}{2009}
\left( \begin{array}{ccc}
          210 & 1645 & 1134 \\
          -1330 & 966 & -1155 \\
          1491 & 630 & -1190 \\
        \end{array}
      \right).
$$

\n Since we have a structure of group on the set of orthogonal
matrices, a natural question is whether one can extend this
algebraic structure to all of the ``equilateral triangles" not only
to those for which $k^2=m^2-mn+n^2$.

\end{document}